\documentclass[10pt]{amsart}
\usepackage{amssymb,amsmath,amsthm,amsfonts,
}

\newcommand{\comment}[1]{}

\newtheorem{lem}{Lemma}
\newtheorem{propn}[lem]{Proposition}

\newtheorem{thm}[lem]{Theorem}

\theoremstyle{remark}

\theoremstyle{definition}
\newtheorem{defn}[lem]{Definition}

\newcommand{\R}{\mathbf R}
\newcommand{\Rn}{{\mathbf R}^d}
\newcommand{\h}{{\mathbf H}}

\DeclareMathOperator{\Rre}{Re}
\newcommand{\p}{\partial}

\newcommand{\vp}{\varphi}

\newcommand{\A}{\alpha}

\newcommand{\B}{\beta}
\newcommand{\lm}{\lambda}

\newcommand{\Tr}{\text{t}}
\DeclareMathOperator{\rank}{rank}
\DeclareMathOperator{\corank}{corank}
\DeclareMathOperator{\Ker}{Ker}

\newcommand{\cc}{\mathcal{C}}
\newcommand{\pr}{\pi_R}
\newcommand{\pl}{\pi_L}

\newcommand{\twxy}{\langle x',y'\rangle}
\newcommand{\twxz}{\langle x',z'\rangle}

\begin{document}

\title[Strongly singular Radon transforms]
{Strongly singular Radon transforms on the Heisenberg group and folding singularities
}
\author{Norberto Laghi\hspace{3cm}Neil Lyall}
\thanks{Both authors were partially supported by HARP grants from the European Commission.}
\address{School of Mathematics, The University of Edinburgh, JCM
 Building, The King's Buildings, Edinburgh EH9 3JZ, United Kingdom}
\email{N.Laghi@ed.ac.uk} 
\address{Department of Mathematics, The University of Georgia, Boyd
  Graduate Studies Research Center, Athens, GA 30602, USA}
\email{lyall@math.uga.edu}
\subjclass{44A12, 42B20, 43A80}
\keywords{Strongly singular integrals, Radon transforms, folding singularities}

\begin{abstract} We prove sharp $L^2$ regularity results
for classes of strongly singular Radon transfoms on the Heisenberg
group by means of oscillatory integrals. We show that the problem in
question can be effectively treated by establishing uniform estimates
for certain oscillatory integrals whose canonical relations project
with two-sided fold singularities; this new approach also allows us to treat
operators which are not necessarily translation invariant.
\end{abstract}
\maketitle

\setlength{\parskip}{5pt}

\section{Introduction}
The principal aim of this work is to study the behaviour of integral operators
acting on functions on the Heisenberg group $\h^n;$ these arise as
natural generalisations of their Euclidean counterparts, often known
as singular Radon transforms. Such integral transforms combine
properties of singular integrals and averages along families of
submanifolds of $\R^d,$ and have attracted great interest in
recent years; for the most recent results and further references see \cite{2}. 

\subsection{Formulation of the problem on the Heisenberg group}

To describe the objects we
shall be interested in, we recall a real-variable characterisation of the
Heisenberg group $\h^n;$ as a topological space this group can be identified
with $\R^{2n+1},$ but Euclidean addition is replaced by the group
operation
\begin{equation}
(x,t)\cdot (y,s)=(x+y,s+t-2\, x^\Tr J y)
\end{equation}
where $J$ denotes the standard symplectic matrix on $\R^{2n}$, namely
\[J=\left(\begin{matrix} 0 &I_n \\ -I_n &0 \end{matrix}\right)\]
and inverses are given by $(x,t)^{-1}=-(x,t)$.
We shall often refer to this last term as simply the twist. The centre of the group $\h^n$ is then given by those elements of the form $(0,t)\in\h^n.$

A problem considered by Geller and Stein in \cite{5} was the following: suppose
$K$ is a Calder\'{o}n-Zygmund kernel in $\R^{2n},$ and $M$ is the
distributional kernel given by the tensor product of $K$ with the Dirac
delta in the central direction, namely 
\[M(x,t)=K(x)\delta(t);\] 
then what are the $L^p$ mapping properties of the singular Radon
transform on $\h^n$ defined by setting $Tf=f*M,$ where convolution is taken with respect to the group structure? Geller and Stein showed that these operators were in fact bounded on  $L^p(\h^n)$ for  $1<p<\infty.$ 

In \cite{8} the second author considered operators $R$ obtained by taking group convolution with the distribution
\begin{equation}\label{dist}
M(x,t)=K_{\A,\B}(x)\delta(t-\phi(x)),
\end{equation}
where $K_{\A,\B}$ is a distribution on $\R^{2n}$ that away from the origin agrees with the function
\begin{equation}\label{K}  
K_{\A,\B}(x)=|x|^{-2n-\A}e^{i|x|^{-\B}}\chi(|x|),
\end{equation}
with $\B>0$ and $\chi$ a smooth cut off function which equals one near the origin\footnote{ \  The distribution-valued function $\A \mapsto K_{\A,\B}$, initially defined for $\Rre \A<0$, continues analytically to the entire complex plane.}. Using group Fourier transform techniques it was shown that if $\phi\equiv0$, or $\phi(x)=|x|^\kappa$ with $\kappa\geq 2$, then $\|R f\|_2\leq C\|f\|_2$ if and only if $\A\leq(n-1/6)\B$. 

Kernels of the form (\ref{K}) were considered by Wainger \cite{10} and C. Fefferman \cite{4} in the context of strongly singular convolution operators; further generalisations
can be found in Lyall \cite{8}. Note that if we choose $\phi\equiv 0,$ then the operators $R$ above are in fact strongly singular analogues of the operator considered by Geller and Stein. 

In this article we shall be principally interested in the study of \emph{strongly singular Radon transforms} (on the Heisenberg group), which we define to be natural generalisations to the non-translation invariant setting of the operators $R$ discussed above as follows; we define these to be operators of the form
\begin{equation}\label{1}Tf(x,t)=
\int_{\R^{2n+1}}K_{\A,\B}(x,y)\left(\int_{\R}e^{i\tau\left[t-s+2\,x^\Tr J y-\phi(x,y)\right]}d\tau\right) f(y,s)dy\,ds,
\end{equation}
where $K_{\A,\B}$ is now a \emph{strongly singular integral kernel}\footnote{ \ Since our operators are not going to be
necessarily translation invariant, the kernel $K_{\A,\B}$ is given by a
distribution on the product of the spaces as defined below.} on $\R^{2n}\times\R^{2n}$. We shall make some specific assumptions on the function $\phi$ later.

Here we shall not aim for the most general
definition of such a kernel; for us 
a strongly singular kernel on $\R^{2n}\times\R^{2n}$
will be a distribution of the form 
\begin{equation}
K_{\A,\B}(x,y)=e^{i |x-y|^{-\beta}}a(x,y)\end{equation} with $\B>0$, where the amplitude
$a$ is supported in a small neighbourhood
of the diagonal $\Delta=\{(x,y)\in\R^{2n}\times\R^{2n}: x=y\},$ is smooth
away from $\Delta$ and satisfies the estimates
\begin{equation}\label{diff}
\left|D^{\mu}_{x,y}a(x,y)\right|\leq C_\mu
|x-y|^{-2n-\alpha-|\mu|}\quad\text{when } x\neq y,
\end{equation}
for every
 multi-index $\mu$; here $\alpha\geq 0$.\footnote{ \ Of course such a
definition, as well as (\ref{K}), is valid also in the odd-dimensional case.}

We shall study (\ref{1}) for two different classes of functions $\phi$ for which we shall make very different qualitative and quantitative assumptions.
Our main result is the following.
\begin{thm}\label{t1} Consider the operator (\ref{1}) with phase function $\phi$ satisfying either of the following conditions:
\begin{itemize}
\item[(i)] $\phi\in C^{\infty}(U\setminus\Delta)$, where $U$ is a neighbourhood
of the diagonal $\Delta\subset\R^{2n}\times\R^{2n}$ with $U\supset supp(a)$, and for some $\kappa>2$ satisfies the differential inequalities
\[\left|D^{\mu}_{x,y}\phi(x,y)\right|\leq C_\mu\left|x-y\right|^{\kappa-|\mu|}\] for all $x\ne y$ and every multiindex $\mu$.
\item[(ii)] $\phi(x,y)=\vp(x-y),$
where $\vp$ is smooth and supported in a small neighbourhood of the origin, with
\[\nabla^2_x\vp(0)=4B\]
where $B=\left(b_i\delta_{i,j}\right)$ with $b_i=b_{i+n}$ a real constant for $i=1,\dots,n$. 
\end{itemize}
Then $T:L^2(\h^n)\to L^2(\h^n)$ if and only if $\alpha\leq (n-1/6)\beta.$
\end{thm}
We note that our second result only concerns operators associated with translation-invariant phase functions. 
The reason for requiring the phase function to have a special form will be clear from the 
arguments provided in the proof.
The model example of such a phase is $\vp(x)=|x|^2$, more generally we can also consider phases of the form $\vp(x)=\sigma(|x|^2),$
where $\sigma$ is a smooth function supported in a neighbourhood of the origin. 
 
We further note that the necessity of the results in Theorem \ref{t1} was shown in \cite{8}.

\subsection{Strongly singular integrals along curves in $\R^d$}

It is standard and well known that the Hilbert transform along curves:
\begin{equation}\label{Hilbert}
H_{\Gamma}f(x)=\text{p.v.}\int_{-1}^1 f(x-\Gamma(t))\frac{dt}{t},
\end{equation}
is bounded on $L^p(\R^d)$, for $1<p<\infty$, where $\Gamma(t)$ is an appropriate curve in $\R^d$. In particular, it was shown by Nagel, Rivi\`ere, and Wainger in \cite{NRW}  that $\|H_{\Gamma}f\|_p\leq C\|f\|_p$, for $1<p<\infty$, where $\Gamma(t)=(t,t|t|^k)$, $k\geq1$, is a curve in $\R^2$, see also Stein and Wainger \cite{StWa}. This work had been originally initiated by Fabes and Rivi\`ere \cite{FR}. 
 
Continuing on the work of Zielinski \cite{Z}, Chandarana \cite{1} studied strongly singular analogues of the above operators, in particular he considered operators on $\R^2$ that take the form
\begin{equation}\label{sh}
Tf(x,t)=\text{p.v.}\int_{-1}^1 H_{\A,\B}(s)f(x-s,t-s|s|^k)ds,\end{equation} where
$H_{\A,\B}(x)=x^{-1}|x|^{-\A}e^{i|x|^{-\beta}}$ is a strongly singular (convolution) kernel in $\R$ which
enjoys some additional cancellation (note that $H_{\A,\B}$ is an odd function for $x\neq 0$).  Note that the convolution kernel $M$ of the operator (\ref{sh}) can of course be written as
\[M(x,t)=H_{\A,\B}(x)\delta(t-x|x|^k)),\] 
which is clearly very reminiscent of (\ref{dist}). 
 
In Section \ref{R2} we shall indicate how the techniques introduced to study operators of the form (\ref{1}) can
be employed to revisit and generalise these results. We however point out that this approach is not exactly necessary and that one can also obtain the result below by simply appealing to van der Corput's lemma, see \cite{curves}.

With our oscillatory integral techniques it is natural to consider operators given by averaging a more general strongly singular kernels over a \emph{smooth} curve $\Gamma(t)=(t,\gamma(t))$. More specifically,  we consider the operators
\begin{equation}\label{euc}
T_\gamma f(x,t)=\int_{\R^2}\int_{\R}
e^{i\left[|x-y|^{-\beta}+\tau(t-s-\gamma(x-y))\right]}a(x,y)d\tau
f(y,s)\,dy\,ds,\end{equation}
where the amplitude $a$ is supported in a small neighbourhood of the diagonal and satisfies the differential inequalities (\ref{diff}) with $n=1/2$. 

\begin{thm}\label{t4} Consider the operator (\ref{euc}) and suppose the smooth curve $\gamma(t)$ has curvature which does not vanish to infinite order in a small neighbourhood of the origin, then
$T_\gamma$ is bounded on $L^2(\R^2)$ if and only if $\alpha\leq \beta/3.$
\end{thm}

\comment{
\begin{thm}\label{t2} Consider  the operator (\ref{1}) and let $\phi(x,y)=\sigma(|x-y|^2),$
where $\sigma$ is smooth in a neighbourhood of the origin in $\R.$ 
Then $T:L^2(\h^n)\to L^2(\h^n)$ for $\alpha\leq(n-1/6)\beta.$
\end{thm}
}

\section{Standard oscillatory integral operator estimates}

Key to our arguments is the following proposition of H\"ormander \cite{7}, \cite{FLp}.
\begin{propn}\label{vcp}
Let $\Psi$ be a smooth function supported on the set $\{(x,y)\in \Rn\times\Rn:|x-y|\leq C\}$ and $\Phi$ be real-valued and smooth on the support of $\Psi$. If we assume that all partial derivatives of $\Psi$ and $\Phi$ are bounded and that
\begin{equation}\label{nondeg}
\det\Bigl(\frac{\p^2\Phi}{\p x_k \p y_\ell}\Bigr)\neq0
\end{equation}
on the support of $\Psi$, then for all $\lm>0$
\[\Bigl\|\int_{\Rn}e^{i\lm\Phi(x,y)}\Psi(x,y)f(y)dy\Bigr\|_{L^2(\Rn)}\leq A(1+\lm)^{-d/2}\|f\|_{L^2(\Rn)}.\]
\end{propn}

Consider the canonical relation  
\[\cc_\Phi=\left\{ \left(x,\Phi_x,y,-\Phi_y\right)\right\}\subset T^*(\R^d_x)\times  T^*(\R^d_y)\]  
associated to the phase function $\Phi$. The non-degeneracy assumption (\ref{nondeg}) is equivalent to the condition that the two projection maps
\[\pl:\cc_\Phi\to T^*(\R^d_x)\quad\text{and}\quad\pr:\cc_\Phi\to T^*(\R^d_y)\]
are local diffeomorphisms.

We also take this opportunity to recall the notion of a map
having fold singularities\footnote{ \ For a detailed and
interesting description of the several kinds of singularities which
are relevant in the theory of oscillatory integral operators
one should consult \cite{3} and \cite{6}.}  and a fundamental result stemming from the work of Melrose and Taylor \cite{9} (see also \cite{PaSo})
which we shall use in this work.
\begin{defn}\label{d1}Let $M_1, M_2$ be smooth manifolds of dimension $n,$ and
let $f:M_1\to M_2$ be a smooth map of $\corank\leq 1.$ Define the
singular variety $S=\{P\in M_1: f\text{ is not locally 1-1 at $P$}\}.$ Then
we say 
that
$f$ has a fold at $P_0$ if
\begin{enumerate}\item[(i)] $\rank\left(Df)\right|_{P_0}=n-1,$ 
\item[(ii)] $\det\,(Df)$ vanishes of first order at $P_0,$
\item[(iii)] $\Ker\left(Df)\right|_{P_0}+T_{P_0}S=T_{P_0}M_1.$
\end{enumerate}
\end{defn}
\begin{propn}[Pan-Sogge]\label{t3} If $\Psi$ and $\Phi$ are, with the exception of Condition (\ref{nondeg}), as in Proposition \ref{vcp}, and $\Phi$ gives rise to a canonical relation whose projections $\pl$ and $\pr$ have at most fold singularities, then for all $\lm>0$
\[\Bigl\|\int_{\Rn}e^{i\lm\Phi(x,y)}\Psi(x,y)f(y)dy\Bigr\|_{L^2(\Rn)}\leq A(1+\lm)^{-d/2+1/6}\|f\|_{L^2(\Rn)}.\] The constant $A$
 depends on the size of the support and the $C^{\infty}$ seminorms of
 $\Psi$, 
as well as the $C^{\infty}$ seminorms of the phase function, remaining bounded if both of these quantities are bounded. The estimates are stable under small perturbations of the phase function in the $C^{\infty}$ topology.
\end{propn}

\section{Decomposition of the operator}\label{decomp} We now introduce decompositions which are convenient in the analysis
of operator (\ref{1}). Let $\zeta$ be a smooth bump function in $C^{\infty}(\R_+)$  with
$\zeta(t)=1\text{ for }t\leq 1/2$ and $\zeta(t)=0\text{ for }t\geq 1$, and define
$\vartheta(t)=\zeta(t)-\zeta(2t)$; then $\sum_{j=1}^{\infty}\vartheta(2^j|t|)\equiv 1$ for $|t|\leq 1/2,\ t\neq 0.$

Next, consider a partition of unity
of the interval $[1/4,1]$ by means of function $\chi_h,$ centred at points $a_h\in [1/4,1]$ with the property that\[\chi_h(t)=\begin{cases}1\text{ if }a_h-\delta\leq t\leq a_h+\delta\\
0\text{ if }a_h-2\delta\leq t\leq a_h+2\delta\end{cases}\]
and
\[\quad\sum_{h=1}^{O(\delta^{-1})}\chi_h(t)=
\begin{cases}1\text{ if }1/4\leq t\leq 1\\ 0\text{ if }1/4-2\delta\leq t\leq 1+2\delta\end{cases}\]
where
$\delta$ is understood to be a small but fixed number.
Note that we have
\[\sum_{h=1}^{O(\delta^{-1})}\chi_h(|t|)\vartheta(|t|)=
\vartheta(|t|).\] 

Further, we decompose the space
$\R^{2n}_y$ into thin half-cones of aperture $\delta$ centred at the point $x$ by means of cutoff functions
$\chi_{\delta}(x,y)$ homogeneous of degree 0;  $O(\delta^{-2n})$ operators are then produced. 

Since both the former and the latter partitions of unity produce a finite number of operators, we shall abuse notation and incorporate the cutoff functions in the amplitude. 

We thus define
\begin{equation}\label{dyadic}T_j f(x,t)=\int_{\R^{2n+1}}\int_{\R} e^{i\left[|x-y|^{-\beta}+\tau(t-s+2\,x^{\Tr}Jy-\phi(x,y))\right]} a_j(x,y)d\tau f(y,s)dy\,ds,\end{equation}
where the amplitude $a_j$ is given by
\[a_j(x,y)=\chi_{\delta}(x,y)\chi_h(2^j|x-y|)a(x,y).\]

\begin{thm}[Key Estimate] \label{Main} If $\phi$ satisfies either Condition (i) or (ii) of Theorem \ref{t1}, then
\[ \|T_{j}f\|_{L^2(\h^n)}\leq C 2^{j(\alpha-(n-1/6)\beta)}\|f\|_{L^2(\h^n)}.\]
\end{thm}   

Theorem \ref{t1} now follows from a standard application of Cotlar's lemma since our operators $T_j$ are, in the following sense, almost orthogonal.

\begin{propn}\label{p2} If $\alpha\leq(n-1/6)\beta$, then the operators $T_{j}$ satisfy the estimate
\[\|T_{j}^*T_{j'}\|_{L^2(\h^n)\to L^2(\h^n)}+\|T_{j'}T_{j}^*\|_{L^2(\h^n)\to L^2(\h^n)}\lesssim 2^{-\beta|j'-j|/6}.\]
\end{propn}

The bulk of the proof of Theorem \ref{Main} is postponed to section \ref{proof}. First we turn our attention to making some additional reductions and establishing Proposition \ref{p2}.

\section{Further reductions and the proof of Proposition \ref{p2}}\label{red}

Taking Fourier transforms in the last variable one obtains the new operator
\[\widetilde{T_j f}(x,\tau)=\int_{\R^{2n}}
 e^{i\left[|x-y|^{-\beta}+\tau(2\,x^{\Tr}Jy-\phi(x,y))\right]} a_j(x,y)\widetilde{f}(y,\tau)dy.\]
It then follows from Plancherel's theorem and rescaling that  establishing Theorem \ref{Main} is equivalent to
verifying that the operators
\begin{equation}\label{sj1}T_{j,\tau}f(x)=2^{j\A}\int_{\R^{2n}} e^{i\left[2^{j\B}|x-y|^{-\beta}+2^{-2j}\tau(2\,x^{\Tr}Jy-2^{2j}\phi(2^{-j}x,2^{-j}y))\right]} b(x,y)f(y)dy\end{equation}
satisfy the estimates
\begin{equation}\label{key}
\|T_{j,\tau}f\|_{L^2(\R^{2n})}\leq C 2^{j(\alpha-(n-1/6)\beta)}\|f\|_{L^2(\R^{2n})}
\end{equation} uniformly in 
$\tau$, where 
\[b(x,y)=2^{-j(2n+\A)}a_j(2^{-j}x,2^{-j}y)\] is smooth, compactly supported and satisfies pointwise estimates which are uniform in $j$.

A further preparatory statement concerns the behaviour of the operator
(\ref{sj1}) when the parameter $2^{-j(\B+2)}|\tau|$ in front of the second term in the phase function is either very large or very small.

\begin{propn} \label{p1} There exists $\epsilon>0$ fixed, such that if 
$2^{-j(\B+2)}|\tau|\notin (\epsilon,\epsilon^{-1})$ then we have
\[\|T_{j,\tau}f\|_{L^2(\R^{2n})}\leq A 2^{j\alpha}\min\{2^{-jn\beta},2^{j2n}|\tau|^{-n}\}\|f\|_{L^2(\R^{2n})}\] with $\epsilon$ and $A$ independent of $j$ and $\tau$. 
\end{propn} 

This result is an immediate consequence of the continuity of the determinant function and Proposition \ref{vcp} once we have established the following two lemmas.

\begin{lem}\label{feff}
Let $\Phi_1(x,y)=|x-y|^{-\B}$, then $\det\bigl(\frac{\p^2\Phi_1}{\p x_k \p y_\ell}\bigr)\neq0$ whenever $x\ne y$ and $\B\neq -1$.
\end{lem}
\begin{proof}
It is easy to verify that 
\[\left(\Phi_1\right)_{xy}(x,y)=\beta |x-y|^{-(\beta+2)}(I-(\beta+2)uu^\Tr),\] 
where $u=(x-y)/{|x-y|}.$ We then employ a device
introduced by C. Fefferman to compute the determinant of this matrix; namely let $R$ be the rotation matrix that takes the vector $u$ to the vector $e_1=(1,0,\ldots,0)\in\R^{2n}.$ Clearly $\det (R)=1$ and we have
\[\det\left(\Phi_1\right)_{xy}(x,y)=\det\left(\beta |x-y|^{-(\beta+2)}(I-(\beta+2)E_{1,1})\right)=-(\B+1)\beta^{2n}|x-y|^{-2n(\beta+2)};\]
here $E_{1,1}$ denotes the matrix whose $(1,1)$ entry is 1, while all the other entries are 0.
\end{proof}
\begin{lem}\label{l1} If $\phi$ satisfies either Condition (i) or (ii) of Theorem \ref{t1} and \[\Phi_2(x,y)=2\,x^{\Tr}Jy-2^{2j}\phi(2^{-j}x,2^{-j}y),\] then $\det\bigl(\frac{\p^2\Phi_2}{\p x_k \p y_\ell}\bigr)\neq0$ whenever $|x-y|\geq c>0$ and $j$ is sufficiently large.
\end{lem}
\begin{proof} It is easy to verify that
\[\left(\Phi_2\right)_{xy}(x,y)=2J-\phi_{xy}(2^{-j}x,2^{-j}y).\] 

If $\phi$ satisfies Condition (i) of Theorem \ref{t1}, then we clearly have that
\[\left(\partial_{x_k y_\ell}\phi\right)(2^{-j}x,2^{-j}y)\leq C 2^{-j(\kappa-2)},\]
for $\kappa>2.$ Consequently the second term is truly an error when $j$ is sufficiently large, and the conclusion follows.

If $\phi$ satisfies Condition (ii) of Theorem \ref{t1}, then it follows from the Taylor expansion
\[\varphi(x)=\varphi(0)+\nabla_{x}\varphi(0)\cdot x+\tfrac{1}{2} x^\Tr \nabla^2_x\vp(0)x+O(|x|^3),\]
that
\[\left(\partial_{x_k y_\ell}\phi\right)(2^{-j}x,2^{-j}y)= -2B+O(2^{-j}).\]
The result then follows in this case from the additional observation that
\[\det\left(2J+2B\right)=\prod_{i=1}^n(4b_i^2+4).\qedhere\]
\end{proof}

We conclude this section by showing that the dyadic operators $T_j$ are almost orthogonal.

\begin{proof}[Proof of Proposition \ref{p2}] We shall only establish the desired estimate for $T^{*}_{j}T_{j'};$
the proof of the other estimate is analogous. 
We again observe that by taking Fourier transforms in the last variables and rescaling it suffices to prove appropriate uniform estimates for the $L^2(\R^{2n})\to L^2(\R^{2n})$ norm of $T_{j,\tau}^*T_{j',\tau}$. 

It follows from Theorem \ref{Main} that the operators $T_{j,\tau}$ are uniformly bounded on $L^2(\R^{2n})$ whenever $\alpha\leq(n-1/6)\beta$, since we also have the trivial estimate
\begin{equation}\label{opest}
\|T_{j,\tau}^*T_{j',\tau}\|\leq\|T_{j,\tau}\| \,\|T_{j',\tau}\|,
\end{equation}
we can clearly assume that $|j'-j|\gg1$.

Let $\epsilon>0$ be the constant given in
Proposition \ref{p1} and without loss in generality we assume that $j'\geq j+C_0,$ where $2^{C_0(\B+2)}\geq\epsilon^{-2}$. 
We now distinguish between two cases.
\begin{itemize}
\item[(i)] If $2^{-j'(\beta+2)}|\tau|\notin [\epsilon,\epsilon^{-1}]$, then it follows from (\ref{opest}) and Proposition \ref{p1} that
\[\|T_{j,\tau}^*T_{j',\tau}\|\leq C\|T_{j',\tau}\|\leq C2^{j'(\alpha-n\beta)}\leq C2^{-j'\B/6}.\]
\item[(ii)] If $2^{-j'(\beta+2)}|\tau|\in [\epsilon,\epsilon^{-1}]$, then $2^{-j(\beta+2)}|\tau|\geq 2^{C_0(\B+2)}\epsilon$, and hence appealing to (\ref{opest}) and Proposition \ref{p1} one more time it follows that
\[\|T_{j,\tau}^*T_{j',\tau}\|\leq C\|T_{j,\tau}\|\leq C2^{j\alpha}2^{j2n}|\tau|^{-n}\leq C2^{-n(j'-j)(\B+2)}.\qedhere\]
\end{itemize}
\end{proof}

\section{Proof of Theorem \ref{Main}}\label{proof}

It follows from the reductions made in Section \ref{red} that in order to prove Theorem \ref{Main} (and hence Theorem \ref{t1}) it suffices to establish estimate (\ref{key}) for the operators $T_{j,\tau}$. 
We recall that
\begin{equation*}T_{j,\tau}f(x)=2^{j\A}\int_{\R^{2n}} e^{i2^{j\B}\left[|x-y|^{-\beta}+2^{-j(\B+2)}\tau\left(2\,x^{\Tr}Jy-2^{2j}\phi(2^{-j}x,2^{-j}y)\right)\right]} b(x,y)f(y)dy\end{equation*}
where $b$ is smooth, compactly supported and satisfies pointwise estimates which are independent of $j$.

We note that if $\phi$ satisfies Condition (i) of Theorem \ref{t1}, then we have that 
\[2^{2j}\phi(2^{-j}x,2^{-j}y)=O(2^{-j(\kappa-2)})\]
where this inequality holds in the $C^m$ topology for any $m\in\mathbf{Z}_{+},$ meaning that the derivatives up to order $m$ also satisfy this bound. 
While if $\phi$ satisfies Condition (ii) of Theorem \ref{t1}, then we may assume to have
\[2^{2j}\phi(2^{-j}x,2^{-j}y)=2\,(x-y)^\Tr B(x-y)+O(2^{-j}),\] 
as in the proof of Lemma \ref{l1}. 
In view of these observation we will first show how the desired bounds are obtained in the case when the errors above are identically zero. 

In light of Proposition \ref{p1} we may assume that 
\[2^{-j(\B+2)}|\tau|\in [\epsilon,\epsilon^{-1}]\] for some $0<\epsilon<1$ fixed. If we assume that $\tau>0$ (the case for $\tau<0$ is similar) and rescale $T_{j,\tau}$ by performing the changes of variables 
\[ x\mapsto2^{j}\tau^{-1/(\B+2)}x,\quad y\mapsto2^{j}\tau^{-1/(\B+2)}y, \]  
we are led, in the case when the errors are identically zero, to study operators of the form\footnote{ \ Note that the 
factors of $2^{j}\tau^{-1/(\B+2)}$ produced by the changes of variables are clearly insignificant and can be neglected.}
\begin{equation}\label{Tlm}
 T_\lm f(x)=\int e^{i \lm\Phi(x,y)}
\Psi(x,y)f(y)dy 
\end{equation}
where $\lm=\tau^{\B/(\B+2)}\sim 2^{j\B}$, \[\Psi(x,y)=b(2^{j}\tau^{-1/(\B+2)}x,2^{j}\tau^{-1/(\B+2)}y)\] and
\begin{align*}
\Phi(x,y)=\begin{cases} |x-y|^{-\beta}+ 2\,x^\Tr J y &\text{ if $\phi$ satisfies Condition (i)} \\
|x-y|^{-\beta}+ 2\,x^\Tr J y-2\,(x-y)^\Tr H(x-y) &\text{ if $\phi$ satisfies Condition (ii)}\end{cases}.\end{align*}

We shall now establish the following result.

\begin{propn}\label{p3} If $T_\lm$ is of the form (\ref{Tlm}) above, then
\[\|T_\lm f\|_{L^2(\R^{2n})}\leq C \lm^{-(n-1/6)}\|f\|_{L^2(\R^{2n}).}\]
\end{propn}   
\begin{proof}
We now consider the canonical relation \[\cc_{\Phi}=\left\{\left(x,\Phi_{x},y,-\Phi_{y}\right)\right\}
\subset T^*(\R^{2n}_x)\times T^*(\R^{2n}_y) \] associated to the operators $T_\lm$, and in particular the two projections
\[ \pl:\cc_{\Phi}\to T^*(\R^{2n}_x), \quad \pr:\cc_{\Phi}\to T^*(\R^{2n}_y) \] to the cotangent
bundles of the base spaces. We wish to show that both projections
$\pl$ and $\pr$ have at most fold singularities as the result then follows from Proposition \ref{t3}, while the estimates may depend on the parameter $2^{j}\tau^{-1/(\B+2)}$, this is no more a matter of concern as this parameter belongs to a bounded set.
We therefore turn our attention to the derivatives $D\pl$ and $D\pr$. These are given by matrices whose determinants coincide (see \cite{7}) and are equal to $\det( \Phi_{xy})(x,y)$.

We shall present here only the arguments in the case where $\phi$ satisfies Condition (ii) of Theorem \ref{t1}, the other case is simpler.
In this case we have 
\[\Phi_{xy}(x,y)=\left(\Phi_1\right)_{xy}(x,y)+2J+2B,\]
where $\Phi_1(x,y)=|x-y|^{-\beta}$. As in the proof of Lemma \ref{feff} we see that
\begin{align*}
\det\left(\Phi_{xy}\right)(x,y)&=\det \left(\beta |x-y|^{-(\beta+2)}(I-(\beta+2)E_{1,1})+2J+2B\right)\\
&=-\left( \B^2(\B+1)Q^2+2\B^2b_1Q-4b_1^2-4\right)\prod_{i=2}^n\left(\left( \B Q+2b_i \right)^2 +4 \right),
\end{align*}
where $Q=|x-y|^{-(\B+2)}$. Note that it is clear from the first equality above that \[\rank\left(\Phi_{xy}\right)\geq 2n-1,\]
thus both
$\pl$ and $\pr$ are maps of corank $\leq 1.$
Furthermore we see that $\det\left(\Phi_{xy}\right)(x,y)$ vanishes if and only if 
\[\B^2(\B+1)Q^2+2\B^2b_1Q=4b_1^2+4.\]
We now consider the variety\footnote{ \ This is clearly diffeomorphic to the singular variety via the parameterization $(x,y)\to(x,\Phi_x,y,-\Phi_y)$; thus in order to study the properties of the singular variety it suffices to study the properties of $\mathfrak{S}$.}
\[\mathfrak{S}=\{(x,y)\in\Psi:\det\left(\Phi_{xy}\right)(x,y)=0\}.\]
It is easy to then verify that
\begin{equation}\label{grad}
\nabla_{x,y}\left.\det\left(\Phi_{xy}\right)\right|_{\mathfrak{S}}=C_\B|x-y|^{-(\beta+3)}\left((\beta+1)|x-y|^{-(\B+2)}+b_1\right)\left(u,-u\right),
\end{equation}
where $C_{\beta}\ne0$ and as in the proof of Lemma \ref{feff} we have set $u=(x-y)/|x-y|$.

 It is now simple to check that $\det\left(\Phi_{xy}\right) \neq 0$ 
whenever  $\nabla_{x,y}\det\left(\Phi_{xy}\right)=0.$ 
Indeed $\nabla_{x,y}\det\left(\Phi_{xy}\right)=0$ if and only if $b_1=-(\beta+1)Q$ which implies
$\left|\det\left(\Phi_{xy}\right)\right|\geq 4^n.$ Thus, the determinant of $\Phi_{xy}$ vanishes of
the first order on $\mathfrak{S}.$

It now only remains for us to verify the third condition contained in Definition \ref{d1}.
We focus our attention on $\pl$, the arguments for $\pr$ are similar.  We now wish to establish the transversality condition 
\begin{equation}\label{cond3}
\Ker\left.(D\pl)\right|_P+T_P S=T_P\cc
\end{equation} 
for $P\in S$; again it will suffice to work with the variety $\mathfrak{S}$.

First we observe that it follows from (\ref{grad}) that the vector $(u,-u)$ is orthogonal to $\mathfrak{S}$ and furthermore note that if $(v,w)=(v_1,\ldots,v_{2n},w_1,\ldots,w_{2n})\in
\Ker(D\pl),$ then necessarily $v=0$. Therefore in order to establish (\ref{cond3}) we need only verify that if $(v,w)\in \Ker(D\pl)$ is nontrivial, then $u\cdot w\neq 0.$

To prove the claim we assume $u\cdot w=0$,
it then follows that if
\begin{equation*} \Phi_{xy}\,w=
\left(\beta |x-y|^{-(\beta+2)}I+2J+2B\right)w=0,
\end{equation*}
then necessarily $w=0$, since \[\det\left(\beta |x-y|^{-(\beta+2)}I+2J+2B\right)\ne0,\] 
a contradiction.
\end{proof}

The complete proof of estimate (\ref{Main}) now also follows, as it is simple to observe that the errors in the phase function, although they may depend on the parameter $2^{j}\tau^{-1/(\B+2)}$, are in fact $O(2^{-j(\kappa -2)})$ and $O(2^{-j})$ respectively in cases (i) and (ii); since $2^{j}\tau^{-1/(\B+2)}$ is bounded and $j$ can be assumed to be large, they can be regarded as small perturbations of the phase function. This shows
\[\|T_{j,\tau}f\|_{L^2(\R^{2n})}\leq C2^{-j(\A-(n-1/6)\beta)}\|f\|_{L^2(\R^{2n})}\]
uniformly in $\tau$, as desired.

\section{Remarks}\label{remarks}
There are a few questions of interest which are not answered in this paper and which we believe deserve further investigation. 

Firstly, it would be of interest to determine an optimal class of smooth functions $\phi$ for which the estimates of Theorem \ref{t1} hold. 
While part (i) of Theorem \ref{t1} is (in our opinion) fairly satisfactory, the results of part (ii)  can possibly be improved; the difficulties are in the calculations needed to understand the behaviour of determinants.

To be more precise, we note that in our arguments in order to compute the determinant of the mixed hessian of the phase function, we need the matrices involved to commute with rotations (or at least with the rotation employed in the proof); while this may not be necessary for the result to hold, it seems like the calculations needed might
be intractable otherwise.

Furthermore, the twist term which is created by group convolution introduces an ``element of curvature" which we wish to preserve; concretely, we wish the matrix $(\Phi_2)_{xy}$ to have maximal rank (the content of Lemma \ref{l1}), a fact used several times in our arguments. This may not be the case if we consider a general, smooth phase.
One should compare this with the Euclidean result of \S \ref{R2} below. 

It would also be of interest to consider a strongly singular kernel with a more general oscillation,
strongly singular integrals with this property are briefly considered in \cite{8}.

\section{Proof of Theorem \ref{t4}}\label{R2}

The necessity of the condition imposed on the indices $\alpha$ and
$\beta$ is essentially in \cite{1} and the sufficiency truly follows the line of our arguments on the Heisenberg group. In order to decompose the operator (\ref{euc}), define cutoff functions $\chi_h,\chi_{\delta}$ and $\vartheta$ as in \S\ref{decomp}, and let
\begin{equation}\label{tj}
T_j f(x,t)=\int_{\R^{2}}\int_{\R} e^{i\left[|x-y|^{-\beta}+\tau(t-s-\gamma(x-y))\right]} a_j(x,y)d\tau f(y,s)dyds,\end{equation}
where $a_j(x,y)=\chi_{\delta}(x,y)\chi_h(2^j|x-y|)a(x,y).$ 

It follows from our assumption that the curve $\gamma$ is not flat that we may assume $\gamma(0)=\gamma'(0)= \cdots =\gamma^{(k-1)}(0)=0$, while $\gamma^{(k)}(0)\ne0$ for some $k\geq2$.

By taking Fourier transforms in the last (second) variable matters again essentially reduce to showing that the (rescaled) operators
\begin{equation}
T_{j,\tau} f(x)=2^{j\A}\int e^{i2^{j\B}\left[|x-y|^{-\beta}-2^{-j(\B+k)}\tau\Phi_3(x-y)\right]} b(x,y)f(y)dy\end{equation}
where $\Phi_3(x)=2^{jk}\gamma(2^{-j}x)$ satisfy the estimates 
\begin{equation}\label{key2}
\|T_{j,\tau}f\|_{L^2(\R)}\leq C 2^{j(\alpha-\beta/3)}\|f\|_{L^2(\R)}
\end{equation} uniformly in 
$\tau$, where 
$b(x,y)=2^{-j(1+\A)}a_j(2^{-j}x,2^{-j}y)$ is smooth, compactly supported and satisfies pointwise estimates which are independent of $j$.

Writing $\Phi_3(x)=\frac{1}{k!}\gamma^{(k)}(0)x^k+O(2^{-j})$ we see that $\Phi_3''(x)\ne0$ on the support of the kernel provided $j$ is large enough. Thus, the analogue of Proposition \ref{p1} follows easily and we may, analogously to our arguments above, assume that the parameter $2^{-j(\B+k)}|\tau|\in[\epsilon,\epsilon^{-1}]$ for some $0<\epsilon<1$ fixed. 

As before we will assume that $\Phi_3(x)=x^k$ and $\tau>0$, the case for $\tau<0$ can again be treated similarly. It then follows, from the uniformity of the estimates of Melrose and Taylor, that matters essentially reduce to establishing that the operators
\begin{equation}\label{Tlm2}
 T_\lm f(x)=\int e^{i \lm\Phi(x-y)}
\Psi(x,y)f(y)dy 
\end{equation}
where $\lm=\tau^{\B/(\B+k)}\sim 2^{j\B}$, $\Psi(x,y)=b(2^{j}\tau^{-1/(\B+k)}x,2^{j}\tau^{-1/(\B+k)}y)$, and $\Phi(x)=|x|^{-\beta}-x^k$ 
give rise to canonical relations which project with at most fold singularities.
But in this setting this is really rather easy and simply amounts to the observation that if $\Phi''(x_0)=0$, then necessarily $\Phi'''(x_0)\ne0$.

This establishes estimate (\ref{key2}) uniformly in $\tau.$ As almost orthogonality also follows as in Proposition \ref{p2}, this concludes the proof.
\bibliographystyle{siam}

\end{document}